\documentclass{amsart}
\usepackage{amssymb}

\newtheorem{thm}{Theorem}[section]
\newtheorem{cor}[thm]{Corollary}
\newtheorem{lemma}[thm]{Lemma}

\theoremstyle{definition}

\theoremstyle{remark}

\numberwithin{equation}{section}

\newcommand{\GS}{Gr\"obner-Shirshov}

\begin{document}

\title{Gr\"obner-Shirshov basis for the braid group in the
Birman-Ko-Lee-Garside generators  }

\author[Bokut]{L. A. Bokut$^*$}

\thanks{$^*$Supported in part by the Russian Fund of
Basic  Research (project 05--01--00230),
the Russian Fund of Leading Scientific Schools,
and the Integration grant of the Siberian Branch of the
RAS N~1.9}

\maketitle

\address{Sobolev Institute of Mathematics, Novosibirsk 630090, Russia}

\email{bokut@math.nsc.ru}

 \begin{abstract}
In this paper, we obtain  Gr\"obner-Shirshov (non-commutative
Gr\"obner) bases for the braid groups in the Birman-Ko-Lee
generators enriched by new ``Garside word" $\delta$ (\cite{BKL98}).
It gives a new algorithm for getting the Birman-Ko-Lee Normal Form
in the braid groups, and thus a new algorithm for solving the word
problem in these groups.
\end{abstract}

\section{Introduction}

In 1998, J.~Birman, K.H.~Ko and S.J.~Lee ~\cite{BKL98} found a new system of generators
of the $n$-string braid group $B_{n}$
and used this system to obtain new solutions of the word problem and the conjugacy problem.
To be more precise, they found a normal form of words
in the braid group
\begin{align*}
&B_{n}=gp\langle a_{ts}, n\geq t>s\geq 1| a_{ts}a_{rq}=a_{rq}a_{ts},
(t-r)(t-q)(s-r)(s-q)>0,\\
&a_{ts}a_{sr}=a_{tr}a_{ts}=a_{sr}a_{tr}, n\geq t>s>r\geq 1\rangle.
\end{align*}

 Recall that the Birman-Ko-Lee generators $a_{ts}$
are the  elements
$$
a_{ts}=(a_{t-1}a_{t-2} \dots a_{s+1})a_s(a_{s+1}^{-1}\dots
a_{t-2}^{-1}a_{t-1}^{-1}),
$$
where $n\geq t>s\geq 1$ and $a_i, n-1\geq i\geq 1 $ are the Artin
generators of $B_n$.

In the paper \cite{BFKS}, a \GS\ basis of the semigroup of positive
braids $B_{n}^+$
in the Artin generators was found.

In the paper, \cite{Bo05}, I found a \GS\ basis of the braid group
$B_{n}$ in the Artin-Garside generators $a_i, 1\leq i \leq n,
\Delta, \Delta^{-1},$ where
$$
\Delta=\Lambda_1\Lambda_2\dots \Lambda_n, \text{ \ with \ }
\Lambda_i=a_i\dots a_1
$$
(see A.F. Garside \cite{Ga69}).

In this paper, I find a \GS\ basis of $B_n$ in the Birman-Ko-Lee
generators enriched by the new "Garside element"
$$
\delta=a_{nn-1}a_{n-1n-2}\dots a_{21}
$$
invented in the same paper ~\cite{BKL98}.
The generators obtained are called the
Birman-Ko-Lee-Garside generators of the braid group.

As a corollary of our result,
the Composition  Lemma (\cite{Sh62},
\cite{Bo72}, \cite{Bo76},  see below)
immediately implies the Birman-Ko-Lee normal form
in the braid group.

\section{Basic notations and results}
 Let us order generators
$$
\delta^{-1}< \delta <a_{ts}<a_{rq}\ \  \text {\ iff \ }\ \  (t,s)<(r,q)
\text { \ lexicographically. \ }
$$

We order words in this alphabet in the deg-lex way  comparing two words
first by theirs degrees
(lengths) and then lexicographically when the degrees are equal.

Instead of $a_{ij}, i>j$, we write simply $(i,j)$ or $(j,i)$. We also set
$$
(t_m,t_{m-1}, \dots ,t_1)= (t_m,t_{m-1})(t_{m-1},t_{m-2})\dots (t_2,
t_1),
$$
where $t_j\neq t_{j+1}, 1\leq j\leq m-1$.

In this notation, the defining relations of $B_n$ can be written as
\begin{align}
&(t_3,t_2,t_1)=(t_2,t_1,t_3)=(t_1,t_3,t_2), t_3>t_2>t_1, \\
&(k,l)(i,j)=(i,j)(k,l), k>l, i>j, k>i,
\end{align}
where either $k>i>j>l$, or $k>l>i>j$.

Let us assume the following notation:
$$
V_{[t_2,t_1]}, \text { where } n\geq t_2> t_1\geq 1,
$$
is a positive word in $(k,l)$ such that $t_2\geq k>l\geq t_1 $.
We can use any capital Latin letter with indices
instead of $V$,
and any appropriate numbers (for example,
$t_3, t_0$ such that $t_3>t_0$) instead of $t_2, t_1$.

We will also use the following notations:
\begin{align*}
&V_{[t_2-1,t_1]}(t_2,t_1)=(t_2,t_1)V_{[t_2-1,t_1]}', t_2>t_1,
\end{align*}
where
$V_{[t_2-1,t_1]}'=(V_{[t_2-1,t_1]})|_{(k,l)\mapsto
(k,l),\text { if } l\neq t_1; (k,t_1)\mapsto (t_2,k)}$;
\begin{align*}
&W_{[t_2,t_1]}(t_1,t_0)=(t_1,t_0)W_{[t_2,t_1]}^\star, t_2>t_1>t_0,
\end{align*}
where
$W_{[t_2,t_1]}^\star=(W_{[t_2,t_1]})|_{(k,l)\mapsto
(k,l),\text { if } l\neq t_1; (k,t_1)\mapsto (k,t_0)}$.

Then
\begin{align*}
&V_{[t_2-1,t_1]}'=(V_{[t_2-1,t_1]}')_{[t_2,t_1+1]},\\
&W_{[t_2,t_1]}^\star=(W_{[t_2,t_1]}^\star)_{[t_2,t_0]}
\end{align*}

\begin{lemma}\label{L00}
Relations
\begin{align*}
&V_{[t_2-1,t_1]}(t_2,t_1)=(t_2,t_1)V_{[t_2-1,t_1]}', t_2>t_1,\\
&V_{[t_2-1,t_1]}(t_1,t_0)=(t_1,t_0)V_{[t_2-1,t_1]}^\star, t_1>t_0,\\
\end{align*}
imply
\begin{align*}
&((V_{[t_2-1,t_1]})^\star)_{[t_2-1,t_0]}(t_2,t_0)=(t_2,t_0)V_{[t_2-1,t_1]}',
\text
{so }((V_{[t_2-1,t_1]}^\star)_{[t_2-1,t_0]})'=V_{[t_2-1,t_1]}';\\
&\text {it means that in the previous notations }\\
&(({V_{[t_2-1,t_1]})^\star})_{[t_2-1,t_0]}'=V_{[t_2-1,t_1]}', \text{
or } {V^\star}'=V' \text { for short }.
\end{align*}
\end{lemma}

\begin{proof}
We have the following letter transformations
$V_{[t_2-1,t_1]}\rightarrow (V_{[t_2-1,t_1]})^\star\rightarrow
((V_{[t_2-1,t_1]})^\star)_{[t_2,t_0]}'$:
$$
(k,l)(t_2-1\geq k>l>t_1)\mapsto(k,l)\mapsto(k,l); (k,t_1)(t_2-1\geq
k>t_1)\mapsto (k,t_0)\mapsto(t_2,k).
$$
\end{proof}

\begin{lemma}\label{L1} Let $t_m>t_{m-1}> \dots > t_1$. Then
$$
(t_m,t_{m-1}, \dots ,t_1)=(t_k, \dots ,t_1, t_m,\dots ,t_{k+1})
$$
for any $1\leq k\leq m$.
\end{lemma}
\begin{proof}For $m=2$, the statement is provided by definition.
Let us proceed by induction on $m\geq 3$.
For $m=3$,
these are defining relations.

$(m-1)\Rightarrow m, m\geq 4$.

Use the induction on $k\geq 1$.

$k=1$:
\begin{align*}
&(t_m,\dots ,t_1)= (t_m,t_{m-1})(t_{m-1},\dots
,t_1)=(t_m,t_{m-1})(t_1,t_{m-1},\dots ,t_2)=\\
&(t_m,t_{m-1})(t_1,t_{m-1})(t_{m-1},\dots
,t_2)=(t_m,t_{m-1},t_1)(t_{m-1},\dots ,t_2)=\\
&(t_1,t_m,t_{m-1})(t_{m-1},\dots ,t_2)=(t_1,t_m,t_{m-1},\dots ,t_2).
\end{align*}

$k\Rightarrow (k+1)$:
\begin{align*}
&(t_m,\dots ,t_1)= (t_k, \dots ,t_1, t_m,\dots ,t_{k+1})=(t_k, \dots
,t_1)(t_1, t_m,\dots
,t_{k+1})=\\
&(t_k, \dots ,t_1)(t_{k+1},t_1, t_m,\dots ,t_{k+2})=
(t_k, \dots ,t_1)(t_{k+1}, t_1)(t_1, t_m,\dots ,t_{k+2})=\\
&(t_k, \dots ,t_1)(t_1,t_{k+1})(t_1, t_m,\dots ,t_{k+2})=
(t_k, \dots ,t_1, t_{k+1})(t_1,t_m, \dots ,t_{k+2})=\\
&(t_{k+1},t_k, \dots ,t_1)(t_1, t_m,\dots ,t_{k+2} )=(t_{k+1},\dots
,t_1, t_m,\dots ,t_{k+2}
).
\end{align*}
\end{proof}

\begin{lemma}\label{L2} Let $t_3>t_2>t_1$. Then
$$
(t_3,t_2, t_1)=(t_2,t_1)(t_3,t_1).
$$
\end{lemma}

\begin{proof}
\begin{align*}
(t_3,t_2,t_1)=(t_2,t_1,t_3)=(t_2,t_1)(t_1,t_3)=(t_2,t_1)(t_3,t_1).
\end{align*}
\end{proof}

\begin{lemma}\label{L3}
Let $t>t_3>t_2>t_1,
t_2>s.$
Then
\begin{align}
&(t,s)(t_2,t_1)(t_3,t_1)=(t_3,t_2)(t,s)(t_2,t_1),\\
&(t_3,s)(t_2,t_1)(t_3,t_1)=(t_2,s)(t_3,s)(t_2,t_1).
\end{align}
\end{lemma}

\begin{proof}
We have
\begin{align*}
&(t,s)(t_2,t_1)(t_3,t_1)=(t,s)(t_3,t_2,t_1)=(t,s)(t_3,t_2)(t_2,t_1)=(t_3,t_2)(t,s)(t_2,t_1),\\
&(t_3,s)(t_2,t_1)(t_3,t_1)=(t_3,s)(t_3,t_2,t_1)=(t_3,s)(t_3,t_2)(t_2,t_1)=(t_2,s)(t_3,s)(t_2,t_1).
\end{align*}
\end{proof}

\begin{lemma}\label{L4} (Lemma 2.3.III, \cite{BKL98}) Let $n\geq t>s\geq 1$. Then
$$
(t,s)\delta=\delta (t+1,s+1),
$$
where $t+1, s+1$ are defined $mod \ n$.
\end{lemma}

\begin{proof} We state the proof to make the exposition complete.
Our proof is different from the one from \cite{BKL98}.
 Also, the cases $n=t$
or $s=1$ in the proof  of Lemma 2.3.III, \cite{BKL98} were given
for readers.

Let $t=n, s=1$:
\begin{align*}
&(n,1)(2,1)\dots (n,1)=(n,1)(n,n-1)(n-1,\dots,
1)=(n-1,1)(n,1)(n-1,\dots, 1)=\\
&(n-1,1)(n-1,\dots, 2)(n,1)(2,1)=(1,n-1,\dots, 2)(n,1)(2,1)=(n-1,\dots,
1)(n,1)(2,1)=\\
&(2,1)\dots (n,1)(2,1).
\end{align*}

Let $t=n$. If $s=n-1$, then
\begin{align*}
&(n,n-1)(2,1)\dots (n-2,1)(n-1,1)(n,1)=(2,1)\dots
(n-2,1)(n,n-1)(n-1,1)(n,1)=\\
&(2,1)\dots (n-2,1)(n-1,1)(n,1)(n,1)=\delta (n,1).
\end{align*}

Take $1<s<n-1$, then
\begin{align*}
&(n,s)(2,1)\dots (s-1,1)(s,1)\dots (n,1)=(2,1)\dots
(s-1,1)(n,s)(s,1)\dots (n,1)=\\
&(2,1)\dots (s-1)(s,1)(n,1)(s+1,1)\dots (n,1)=\\
&(2,1)\dots (s,1)(n,1)(s+1,1)\dots (n-1,1)(n,1)=\\
&(2,1)\dots (s,1)(n,1)(s+1,1)\dots (n-2,1)(n,n-1)(n-1,1)=\\
&(2,1)\dots (s,1)(n,1)(n,n-1)(s+1,1)\dots (n-1,1)=\\
&(2,1)\dots (s,1)(n,1)(n-1,n-2,\dots ,s+1,1)=\\
&(2,1)\dots (s,1)(n-1,1)(n,1)(s+1,1)\dots (n-1,1)=\\
&(2,1)\dots (s,1)(n-1,1)(n,1)(n-1,\dots ,s+1,1)=\\
 &(2,1)\dots
(s,1)(n-1,1)(n-1,n-2)(n-2,n-3)\dots
(s+2,s+1)(n,1)(s+1,1)=\\
&(2,1)\dots (s,1)(s+1,1)\dots (n-1,1)(n,1)(s+1,1)=\delta (s+1,1),
\end{align*}
where we use $(n-1,1)(n-1,n-2)(n-2,n-3)\dots (s+2,s+1)=(1,n-1,\dots,
s+1)=$
$(n-1,\dots, s+1,1)=(s+1,1)\dots (n-1,1).$

Now let $t<n$.

$s=1$:
\begin{align*}
&(t,1)\delta=(t,1)(n,\dots ,t+1,t,\dots ,1)=(t,1)(t,\dots
,1)(1,n,\dots ,t+1)=\\
&(t,1)(t,\dots, 2)(2,1)(n,\dots,t+1, 1)=(1,t,\dots,
2)(2,1)(n,\dots,t+1, 1)=\\
&(t,\dots, 1)(n,\dots, t+1)(2,1)(t+1,1)=(t,\dots, 1)(n,\dots,
t+1)(t+1,1)(t+1,2)=\\
&(t,\dots, 1)(n,\dots, t+1,1)(t+1,2)=(2,1)\dots (t,1)(t+1,1)\dots
(n,1)(t+1,2).
\end{align*}

For $n>t>s>1$ we have:

\begin{align*}
&(t,s)\delta= (t,s)(2,1)\dots (s-1,1)(s,1)\dots (n,1)=(2,1)\dots
(s-1,1)(t,s)(s,1)\dots (n,1)=\\
&(2,1)\dots (s,1)(t,1)(s+1,1)\dots (t,1)(t+1,1)\dots (n,1).\\
&\delta (t+1,s+1)=(2,1)\dots (s,1)(s+1,1)\dots
(t+1,1)(t+1,s+1)(t+2,1)\dots (n,1)=\\
&(2,1)\dots (s,1)(s+1,1)\dots (t,1)(s+1,1)(t+1,1)\dots (n,1).
\end{align*}

So we need to prove that
$$
(t,1)(s+1,1)\dots (t,1)=(s+1,1)\dots (t,1)(s+1,1).
$$

We have
\begin{align*}
&(t,1)(s+1,1)\dots (t,1)=(t,1)(t,t-1,\dots,s+1,1)=\\
&(t,1)(t,t-1)(t-1,\dots,s+1,1)=
(t-1,1)(t,1)(t-1,\dots, s+1,1)=\\
&(t-1,1)(t-1,\dots, s+1)(t,1)(s+1,1)=\\
&(1,t-1,\dots, s+1)(t,1)(s+1,1).\\
&(s+1,1)\dots (t,1)(s+1,1)=(t,t-1,\dots,s+1,1)(s+1,1)=\\
&(t-1,\dots,s+1,1,t)(s+1,1)=(t-1,\dots,s+1,1)(t,1)(s+1,1).
\end{align*}

We are done.
\end{proof}

\begin{lemma}\label{L5}
\begin{align}
&(2,1)V_2(3,1)\dots V_{n-1}(n,1)=\delta V_1'\dots V_{n-1}'
         \label{E5'},
\end{align}
where $V_i=V_{i[i-1,1]}, 2\leq i\leq (n-1)$.
\end{lemma}
\begin{proof}

\end{proof}

Through the paper, we fix notations
$$
n\geq t_3>t_2>t_1\geq 1.
$$

\begin{thm}
A \GS\ basis of $B_n$ in the Birman-Ko-Lee-Garside generators consists
of the
following relations:

\begin{align}
&(k,l)(i,j)=(i,j)(k,l), k>l>i>j,\label{E1}\\
&(k,l)V_{[j-1,1]}(i,j)=(i,j)(k,l)V_{[j-1,1]}, k>i>j>l,\label{E2}\\
&(t_3,t_2)(t_2,t_1)=(t_2,t_1)(t_3,t_1),\label{E3}\\
&(t_3,t_1)V_{[t_2-1,1]}(t_3,t_2)=(t_2,t_1)(t_3,t_1)V_{[t_2-1,1]},\label{E4}\\
&(t,s)V_{[t_2-1,1]}(t_2,t_1)W_{[t_3-1,t_1]}(t_3,t_1)=(t_3,t_2)(t,s)V_{[t_2-1,1]}(t_2,t_1)W_{[t_3-1,t_1]}', \label{E5}\\
&(t_3,s)V_{[t_2-1,1]}(t_2,t_1)W_{[t_3-1,t_1]}(t_3,t_1)=(t_2,s)(t_3,s)V_{[t_2-1,1]}(t_2,t_1)W_{[t_3-1,t_1]}',\label{E6}\\
&(2,1)V_{2[2,1]}(3,1)\dots V_{n-1[n-1,1]}(n,1)=\delta V_{2[2,1]}'\dots V_{n-1[n-1,1]}',\label{E7}\\
&(t,s)\delta=\delta(t+1,s+1),(t,s)\delta^{-1}=\delta^{-1}(t-1,s-1),\label{E8}, t\pm1, s\pm1 (mod \ n)\\
&\delta \delta^{-1}=1, \delta^{-1}\delta =1,\label{E9}
\end{align}
where $V_{[k,l]}$ means as before any word in $(i,j)$ such that
$k\geq i>j\geq l, t>t_3, t_2>s$.
\end{thm}

Recall that a subset $S$ of the free algebra $k\langle X \rangle$
over a field $k$ on $X$ is called a \GS\ set (basis) if every
composition of elements of $S$ is trivial. This definition goes back
to  Shirshov's 1962 paper \cite{Sh62}. We recall the definition of
triviality of a composition below.

For $n\geq t>s\geq 1$ we have (see \cite{BKL98})

$$
(n,\dots,t,s-1,\dots, 1)(t-1,\dots,s)(t,s)= \delta.
$$

Indeed, $(t-1,\dots,s)(t,s)=(t-1,\dots,s,t)=(t,\dots,s)$,
$(n,\dots,t,s-1,\dots, 1)=(s-1,\dots, 1,n,\dots,t)$, $(s-1,\dots,
1,n,\dots,t)(t,\dots,s)=\delta$
(by Lemma \ref{L1}).

As the result, we do not need the letters $(t,s)^{-1}$ and the
relations $(t,s)(t,s)^{-1}=1$, $(t,s)^{-1}(t,s)=1$ in the above
presentation
of the group $B_n$.

\section{Proof of the Theorem}

It is easy to see that formulas (\ref{E1})-(\ref{E9}) are valid in
$B_n$.

We need to prove that all compositions of relations
(\ref{E1})--(\ref{E9}) are trivial.

By "a word" we will mean a positive word in $(i,j), \delta$;
 $u=v$ is either the equality in $B_n^+$ or the graphical equality
(the meaning would be clear from the context).

 We use the following notation for words $u, v$:
$$
u\equiv v,
$$
if $u$ can be transformed to $v$ by the eliminations of leading
words of relations (\ref{E1})--(\ref{E9}), i.e., by the eliminations of
left parts of these relations. Actually, we will use an expansion of
this notation meaning that $u\equiv v$ if
$$
u\mapsto u_1\mapsto u_2\mapsto \dots \mapsto u_k= v,
$$
where $u_i< u$ for all $i$ and each transformation is an application
of (\ref{E1})--(\ref{E9}) (so, in general, only the first
transformation
$u\mapsto u_1$ is the elimination
of the leading word of (\ref{E1})--(\ref{E9})).

Another expansion of that formula is
$$u\equiv v(mod\ w)$$
 meaning that $u$ can be transformed to $v$ as
before and all $u_i<w, u<w$.

By abuse of notations, we
assume that in a word equivalence chain starting with a word $u$,
$$
u\equiv v\equiv w\equiv t \dots
$$
each equivalence $v\equiv w, \ w\equiv t, \dots $ is $mod\ u$.

This agrees with the definition of
triviality of a composition (see \cite{Bo72}, \cite{Bo76}). Namely,
a composition $(f,g)_w$ is called trivial $mod (S,w) $,  if
$$
(f,g)_w=\sum \alpha_ia_is_ib_i, a_i\overline{s_i}b_i<w,
$$
where $s_i\in S, a_i,b_i\in X^*, \alpha_i\in k $.
Here $k\langle X
\rangle$ is a free associative algebra over a field $k$ on a set
$X$, $S\subset k\langle X \rangle$, $X^*$ is the set of all words in
$X$, $\overline{s}$ is the leading monomial of a polynomial $s$.
Recall that
$$
(f,g)_w=fb-ag, w=\overline{f}b=a\overline{g}, deg(f)+deg(g)>deg (w),
$$
or
$$
(f,g)_w=f-agb, w=\overline{f}=a\overline{g}b.
$$
are called the compositions of intersection and including
respectively. The first composition we denote also $f\wedge g$, the
second -- $f\vee g$.

Here $w$ is called the ambiguity of the composition $(f,g)_w$, $a,
b\in X^*$.

Let $S$ be the set of polynomial corresponding to semigroup
relations $u_i=v_i, u_i>v_i $, $(f,g)_w=u-v$ is a composition, $
u,v\in X^*, u,v<w$. The triviality  of $(f,g)_w$ $mod(S,w)$ means
that in the previous sense
$$
u\equiv t (mod\ w), v\equiv t (mod\ w)
$$
for some word $t$.
 Now we need to prove that all compositions of relations
 (\ref{E1})-(\ref{E9}) are trivial.

 Recall that we fix order of any integers $t_3,t_2,t_1$,
$$
  n\geq t_3>t_2>t_1\geq 1.
  $$
We will assume also that
$$
t_4>t_3, t_1>t_0, s_3>s_2>s_1.
$$
Any composition $(f,g)_w$ has a form
$$
(f,g)_w=-(u-v),
$$
where $fb=w-u, ag=w-v$ in the case of composition of intersection,
and $f=w-u, agb=w-v$ in the case of composition of including. We
will use this notation freely.

 Let us consider compositions of (\ref{E1}) with all
others relations. We start with listening all intersection
ambiguities of (\ref{E1}):
\begin{align*}
&(\ref{E1})\wedge(\ref{E1}) (k,l)(i,j)(i_1j_1), k>l>i>j>i_1>j_1,\\
&(\ref{E1})\wedge(\ref{E2}) (t,s)(k,l)V_{[j-1]}(i,j), t>s>k>i>j>l,\\
&(\ref{E2})\wedge(\ref{E1})(k,l)V_{[j-1,1]}(i,j)(i_1,j_1), k>i>j>l>i_1>j_1,\\
&(\ref{E1})\wedge(\ref{E3}) (k,l)(t_2,t_1)(t_2,t_1), k>l>t_3,\\
&(\ref{E3})\wedge(\ref{E1}) (t_3,t_2)(t_2,t_1)(i,j), t_1>i>j,\\
&(\ref{E1})\wedge(\ref{E4}) (k,l)(t_3,t_1)V_{[t_2-1]}(t_3,t_2), k>l>t_3,\\
&(\ref{E4})\wedge(\ref{E1}) (t_3,t_1)V_{[t_2-1]}(t_3,t_2)(i,j), t_1>i>j,\\
&(\ref{E1})\wedge(\ref{E5}) (k,l)(t,s)V_{[t_2-1,1]}(t_2,t_1)W_{[t_3-1,t_1]}(t_3,t_1),k>l>t_3\\
&(\ref{E5})\wedge(\ref{E1}) (t,s)V_{[t_2-1,1]}(t_2,t_1)W_{[t_3-1,t_1]}(t_3,t_1)(i,j), t>s>t_3,t_1>i>j,\\
&(\ref{E1})\wedge(\ref{E6}) (k,l)(t_3,s)V{[t_2-1,1]}(t_3,t_1)(W_{[t_3-1,t_1]}(t_3,t_1),k>l>t_3,\\
&(\ref{E6})\wedge(\ref{E1}) (t_3,s)V_{[t_2-1,1]}(t_3,t_1)W_{[t_3-1,t_1]}(t_3,t_1)(i,j), t_1>i>j, t_2>s,\\
&(\ref{E1})\wedge(\ref{E7}) (k,l)(2,1)V_{1[2,1]}(3,1)\dots V_{n-2[n-1,1]}(n,1), k>l>2\\
&(\ref{E1})\wedge(\ref{E8})(k,l)(i,j)\delta^{\pm 1}, k>l>i>j.
\end{align*}

Let us check three of these compositions as examples:
\begin{align*}
& (\ref{E2})\wedge(\ref{E1}) w=(k,l)V_{[j-1,1]}(i,j)(i_1,j_1), k>i>j>l>i_1>j_1,\\
&u-v=(i,j)(k,l)V_{[j-1,1]}(i_1,j_1)-(k,l)V_{[j-1,1]}(i_1,j_1)(i,j),\\
&V_{[j-1,1]}(i_1,j_1)=V_{1[j-1,1]}, v\equiv
(i,j)(k,l)V_{[j-1,1]}(i_1,j_1);\\
&(\ref{E5})\wedge(\ref{E1}) (t,s)V_{[t_2-1,1]}(t_2,t_1)W_{[t_3-1,t_1]}(t_3,t_1)(i,j), t>s>t_3,t_1>i>j,\\
&u-v=(t_3,t_2)(t,s)V_{[t_2-1,1]}(t_2,t_1)W_{[t_3-1,t_1]}'(i,j)-\\
&(t,s)V_{[t_2-1,1]}(t_2,t_1)W_{[t_3-1,t_1]}(i,j)(t_3,t_1),\\
&v\equiv(t,s)(t_3,t_2)V_{[t_2-1,1]}(t_2,t_1)W_{[t_3-1,t_1]}'(i,j)\equiv\\
&(t_3,t_2)(t,s)V_{[t_2-1,1]}(t_2,t_1)W_{[t_3-1,t_1]}'(i,j);\\
&(\ref{E1})\wedge(\ref{E8})w=(k,l)(i,j)\delta, k>l>i>j,\\
&u-v=(i,j)(k,l)\delta-(k,l)\delta(t+1,s+1),
u\equiv(t,s)\delta(k+1,l+1)\equiv \\
&\delta(t+1,s+1)(k+1,l+1),v\equiv \delta(k+1,l+1)(t+1,s+1)\equiv\\
&\delta(t+1,s+1)(k+1,l+1).
\end{align*}
Here we use that $(k+1,l+1)(i+1,j+1)=(i+1,j+1)(k+1,l+1), n\geq
k>l>i>j\geq 1$ (case $k=n$ should be treated separately).

We proceed with intersection compositions of (\ref{E2}) with
(\ref{E2}), \dots, (\ref{E9}). The ambiguities are:
\begin{align*}
&(\ref{E2})\wedge(\ref{E2})(k,l)V_{[j-1,1]}(i,j)W_{[s-1,1]}(t,s), k>i>t>s>j>l, \\
&(\ref{E2})\wedge(\ref{E3})(k,l)V_{[t_2-1,1]}(t_,t_2)(t_2,t_1), k>t_3, t_2>l, \\
&(\ref{E3})\wedge(\ref{E2}) (t_3,t_2)(t_2,t_1)V_{[j-1,1]}(i,j), t_2>i>j>t_1,\\
&(\ref{E2})\wedge(\ref{E4}) (k,l)V_{[t_1-1,1]}(t_3,t_1)W_{[t_2-1,1]}(t_3,t_2),k>t_3,t_1>l,\\
&(\ref{E4})\wedge(\ref{E2})(t_3,t_1)V_{[t_2-1]}(t_3,t_2)W_{[j-1,1]}(i,j), t_3>i>j>t_2,\\
&(\ref{E2})\wedge(\ref{E5}) (k,l)V_{[s-1,1]}(t,s)W_{2-1,1]}(t_2,t_1)W_{[t_3-1,t_1]},k>t>s>l, t_2>s,\\
&(\ref{E5})\wedge(\ref{E2})(t,s)V_{[t_2-1,1]}(t_2,t_1)W_{[t_3-1,t_1]}(t_3,t_1)R_{[j-1,1]}(i,j),t>t_3>i>j>, t_2>s\\
&(\ref{E2})\wedge(\ref{E6})(k,l)V_{[s-1,1]}(t_3,s)W_{[t_2-1,1]}(t_2,t_1)R_{[t_3-1,t_1]}(t_3,t_1),k>t_3, t_2>l, t_2>s,\\
&(\ref{E6})\wedge(\ref{E2})(t_3,s)V_{[t_2-1,1]}(t_2,t_1)W_{[t_3-1,t_1]}(t_3,t_1)R_{[j-1,1]}(i,j), t_3>i>j>t_1, t_2>s,\\
&(\ref{E7})\wedge(\ref{E2})(2,1)V_{2[2,1]}(3,1)\dots V_{n-1[n-1,1]}(n,1)W_{[j-1,1]}(i,j),n>i>j>1,\\
&(\ref{E2})\wedge(\ref{E6})(k,l)V_{[j-1,1]}(i,j)\delta^{\pm1},k>i>j>l.
 \end{align*}

Let us check two of these compositions as examples:
\begin{align*}
&(\ref{E2})\wedge(\ref{E4})
w=(k,l)V_{[t_1-1,1]}(t_3,t_1)W_{[t_2-1,1]}(t_3,t_2),k>t_3,t_1>l,\\
&u-v=(t_3,t_1)(k,l)V_{[t_1-1,1]}W_{[t_2-1,1]}(t_3,t_2)-(k,l)V_{[t_1-1,1]}(t_2,t_1)(t_3,t_1)W_{[t_2-1,1]},\\
&u\equiv(t_3,t_1)(t_3,t_2)(k,l)V_{[t_1-1,1]}W_{[t_2-1,1]}\equiv(t_2,t_1)(t_3,t_1)(k,l)V_{[t_1-1,1]}W_{[t_2-1,1]},\\
&v\equiv(t_2,t_1)(k,l)V_{[t_1-1,1]}(t_3,t_1)W_{[t_2-1,1]}\equiv(t_2,t_1)(t_3,t_1)(k,l)V_{[t_1-1,1]}W_{[t_2-1,1]};\\
&(\ref{E6})\wedge(\ref{E2})w=(t_3,s)V_{[t_2-1,1]}(t_2,t_1)W_{[t_3-1,t_1]}(t_3,t_1)V_{[j-1,1]}(i,j),\\
&t_3>i>j>t_1, t_2>s,\\
&u-v=(t_2,s)(t_3,s)V_{[t_2-1,1]}(t_2,t_1)W_{[t_3-1,t_1]}'V_{[j-1,1]}(i,j),\\
&v\equiv(t_3,s)V_{[t_2-1,1]}(t_2,t_1)W_{[t_3-1,t_1]}(i,j)(t_3,t_1)V_{[j-1,1]}\equiv\\
&(t_2,s)(t_3,s)V_{[t_2-1,1]}(t_2,t_1)W_{[t_3-1,t_1]}'(i,j)V_{[j-1,1]}\equiv\\
&(t_2,s)(t_3,s)V_{[t_2-1,1]}(t_2,t_1)W_{[t_3-1,t_1]}'V_{[j-1,1]}(i,j).
\end{align*}

We proceed with intersection compositions of (\ref{E3}) with
(\ref{E3}), \dots, (\ref{E9}). The ambiguities are:
\begin{align*}
&(\ref{E3})\wedge(\ref{E3})(t_4,t_3)(t_3,t_2)(t_2,t_1),  \\
&(\ref{E3})\wedge(\ref{E4})(t_4,t_3)(t_3,t_1)V_{[t_2-1,1]}(t_3,t_1),  \\
&(\ref{E4})\wedge(\ref{E3}) (t_3,t_1)V_{[t_2-1,1]}(t_3,t_2)(t_2,t_1), t_2>t_1,\\
&(\ref{E3})\wedge(\ref{E5}) (k,t)(t,s)W_{[t_2-1,1]}(t_2,t_1)W_{[t_3-1,t_1]},k>t,\\
&(\ref{E5})\wedge(\ref{E3})(t,s)V_{[t_2-1,1]}(t_2,t_1)W_{[t_3-1,t_1]}(t_1,t_0),\\
&(\ref{E3})\wedge(\ref{E6}) (t_4,t_3)(t_3,s)W_{[t_2-1,1]}(t_2,t_1)R_{[t_3-1,t_1]}(t_3,t_1)\\
&(\ref{E6})\wedge(\ref{E3})(t_3,s)V_{[t_2-1,1]}(t_2,t_1)W_{[t_3-1,t_1]}(t_3,t_1)(t_1,t_0),t_2>s,\\
&(\ref{E3})\wedge(\ref{E7})(t,2)(2,1)V_{2[2,1]}(3,1)\dots V_{n-1[n-1,1]}(n,1),t>2\\
&(\ref{E3})\wedge(\ref{E8})(t_3,t_2)(t_2,t_1)(i,j)\delta^{\pm1},k>i>j>l.
 \end{align*}

Let us check two of these compositions as examples.

We choose the following compositions; here we use Lemma \ref{L00}:
\begin{align*}
&(\ref{E5})\wedge(\ref{E3})w=(t,s)V_{[t_2-1,1]}(t_2,t_1)W_{[t_3-1,t_1]}(t_3,t_1)(t_1,t_0),\\
&u-v=(t_3,t_2)(t,s)V_{[t_2-1,1]}(t_2,t_1)W_{[t_3-1,t_1]}'(t_1,t_0)-\\
&(t,s)V_{[t_2-1,1]}(t_2,t_1)W_{[t_3-1,t_1]}(t_1,t_0)(t_3,t_0),\\
&u\equiv(t_3,t_2)(t,s)V_{t_2-1,1]}(t_2,t_1)W_{[t_3-1,t_1]}'(t_1,t_0)\equiv\\
&(t_3,t_2)(t,s)V_{[t_2-1,1]}(t_1,t_0)(t_2,t_0)W_{[t_3-1,t_1]}',\\
&v\equiv(t,s)V_{[t_2-1,1]}(t_2,t_1)(t_1,t_0)((W_{[t_3-1,t_1]})^\star)_{[t_3-1,t_0]}(t_3,t_0)\equiv\\
&(t,s)V_{[t_2-1,1]}(t_1,t_0)(t_2,t_0)((W_{[t_3-1,t_1]})^\star)_{[t_3-1,t_0]}(t_3,t_0)\equiv\\
&(t_3,t_2)(t,s)V_{[t_2-1,1]}(t_1,t_0)(t_2,t_0)((W_{[t_3-1,t_1]})^\star)_{[t_3-1,t_0]}'\equiv\\
&(t_3,t_2)(t,s)V_{[t_2-1,1]}(t_1,t_0)(t_2,t_0)W_{[t_3-1,t_1]}';\\
&(\ref{E6})\wedge(\ref{E3})(t_3,s)V_{[t_2-1,1]}(t_2,t_1)W_{[t_3-1,t_1]}(t_3,t_1)(t_1,t_0),t_3>s,\\
&u-v=(t_2,s)(t_3,s)V_{[t_2-1,1]}(t_2,t_1)W_{[t_3-1,t_1]}'(t_1,t_0)-\\
&(t_3,s)V_{[t_2-1,1]}(t_2,t_1)W_{[t_3-1,t_1]}(t_1,t_0)(t_3,t_0),\\
&u\equiv(t_2,s)(t_3,s)V_{[t_2-1,1]}(t_2,t_1)(t_1,t_0)W_{[t_3-1,t_1]}'\equiv\\
&(t_2,s)(t_3,s)V_{[t_2-1,1]}(t_1,t_0)(t_2,t_0)W_{[t_3-1,t_1]}',\\
&v\equiv(t_3,s)V_{[t_2-1,1]}(t_2,t_1)(t_1,t_0)W_{[t_3-1,t_1]}^\star(t_3,t_0)\equiv\\
&(t_3,s)V_{[t_2-1,1]}(t_1,t_0)(t_2,t_0)(W_{[t_3-1,t_1]})^\star(t_3,t_0)\equiv\\
&(t_2,s)(t_3,s)V_{[t_2-1,1]}(t_1,t_0)(t_2,t_0)((W_{[t_3-1,t_1]})^\star)'\equiv\\
&(t_2,s)(t_3,s)V_{[t_2-1,1]}(t_1,t_0)(t_2,t_0)W_{[t_3-1,t_1]}'.\\
\end{align*}

Our next compositions will be (\ref{E4}) with (\ref{E4})-(\ref{E9}).
The ambiguities of intersection are the following:

\begin{align*}
&(\ref{E4})\wedge(\ref{E4})(t_4,t_1)V_{[t_2-1,1]}(t_4,t_2)W_{[t_3-1,1]}(t_4,t_3),\\
&(\ref{E4})\wedge(\ref{E5})(t,s')V_{[s-1,1]}(t,s)W_{[t_2-1,1]}(t_2,t_1)R_{[t_3-1,t_1]}(t_3,t_1),t>t_3, t_2>s>s',\\
 &(\ref{E5})\wedge(\ref{E4})(t,s)V_{[t_2-1,1]}(t_2,t_1)W_{[t_3-1,t_1]}(t_3,t_1)R_{[t_3'-1,1]}(t_3,t_3'),t_3>t_3'>t_1,\\
 &(\ref{E4})\wedge(\ref{E6})(t_3,s')V_{[s-1,1]}(t_3,s)W_{[t_2-1,1]}(t_2,t_1)R_{[t_3-1,t_1]}(t_3,t_1),t_2>s>s',\\
 &(\ref{E6})\wedge(\ref{E4})(t_3,s)V_{[t_2-1,1]}(t_2,t_1)W_{[t_3-1,t_1]}(t_3,t_1)R_{[t_3'-1,1]}(t_3,t_3'), t_3>t_3'>t_1,\\
&(\ref{E7})\wedge(\ref{E4})(2,1)V_{2[2,1]}(3,1)\dots V_{n-1[n-1,1]}(n,1)V_{[k-1,1]}(n,k), n>k>1,\\
&(\ref{E4})\wedge(\ref{E8})(t_3,t_1)V_{[t_2-1,1]}(t_3,t_2)\delta^{\pm1}.
\end{align*}

Let us check one of the compositions:

\begin{align*}
&(\ref{E5})\wedge(\ref{E4})w=(t,s)V_{[t_2-1,1]}(t_2,t_1)W_{[t_3-1,t_1]}(t_3,t_1)R_{[t_3'-1,1]},\\
&(t_3,t_3'),t_3>t_3'>t_1,\\
&u-v=(t_3,t_2)(t,s)V_{[t_2-1,1]}(t_2,t_1)W_{[t_3-1,t_1]}'R_{[t_3'-1,1]}(t_3,t_3')-\\
&(t,s)V_{[t_2-1,1]}(t_2,t_1)W_{[t_3-1,t_1]}(t_3',t_1)(t_3,t_1)R_{[t_3'-1,1]},\\
&v\equiv(t_3,t_2)(t,s)V_{[t_2-1,1]}(t_2,t_1)W_{[t_3-1,t_1]}'(t_3,t_3')R_{[t_3'-1,1]}\equiv\\
&(t_3,t_2)(t,s)V_{[t_2-1,1]}(t_2,t_1)W_{[t_3-1,t_1]}'R_{[t_3'-1,1](t_3,t_3')}.
\end{align*}

Now we proceed with intersection compositions of (\ref{E5}) with
(\ref{E5}) -(\ref{E9}). The ambiguities are the following:

\begin{align*}
&(\ref{E5})\wedge(\ref{E5})(t,s)V_{[t_2-1,1]}(t_2,t_1)W_{[t_3-1,t_1]}(t_3,t_1)\times\\
&R_{[s_2-1,1]}(s_2,s_1)T_{[s_3-1,s_1]}(s_3,s_1),t>t_3, t_2>s, t_3>s_3, s_2>t_1,\\
&(\ref{E5})\wedge(\ref{E6})(t,s)V_{[t_2-1,1]}(t_2,t_1)W_{[t_3-1,t_1]}(t_3,t_1)\times \\
&R_{[s_2-1,1]}(s_2,s_1)T_{[t_3-1,s_1]}(t_3,s_1),t>t_3, t_2>s, t_3>s_2>t_1,\\
&(\ref{E6})\wedge(\ref{E5})(t_3,s)V_{[t_2-1,1]}(t_2,t_1)W_{[t_3-1,t_1]}(t_3,t_1)\times\\
&R_{[s_2-1,1]}(s_2,s_1)T_{[s_3-1,s_1]}(s_3,s_1), t_3>s_3,
s_2>t_1,\\
&(\ref{E7})\wedge(\ref{E5})(2,1)V_{2[2,1]}(3,1)\dots V_{n-1[n-1,1]}(n,1)\times\\
&V_{[t_2-1,1]}(t_2,t_1)W_{[t_3-1,t_1]}(t_3,t_1), t>t_3,\\
&(\ref{E5})\wedge(\ref{E8})(t,s)V_{[t_2-1,1]}(t_2,t_1)W_{[t_3-1,t_1]}(t_3,t_1)\delta^{\pm1}.
\end{align*}

Let us check one composition:

\begin{align*}
&(\ref{E5})\wedge(\ref{E5})w=(t,s)V_{[t_2-1,1]}(t_2,t_1)W_{[t_3-1,t_1]}(t_3,t_1)\times\\
&R_{[s_2-1,1]}(s_2,s_1)T_{[s_3-1,s_1]}(s_3,s_1),t>t_3, t_2>s, t_3>s_3, s_2>t_1,\\
&u-v=(t_3,t_2)(t,s)V_{[t_2-1,1]}(t_2,t_1)W_{[t_3-1,t_1]}'\times\\
&R_{[s_2-1,1]}(s_2,s_1)T_{[s_3-1,s_1]}(s_3,s_1)-(t,s)V_{[t_2-1,1]}(t_2,t_1)W_{[t_3-1,t_1]}(s_3,s_2)(t_3,t_1)\times\\
&R_{[s_2-1,1]}(s_2,s_1)T_{[s_3-1,s_1]}',v\equiv(t_3,t_2)(t,s)V_{[t_2-1,1]}(t_2,t_1)W_{[t_3-1,t_1]}'(s_3,s_2)\times\\
&R_{[s_2-1,1]}(s_2,s_1)T_{[s_3-1,s_1]}'\equiv(t_3,t_2)(t,s)V_{[t_2-1,1]}(t_2,t_1)W_{[t_3-1,t_1]}'R_{[s_2-1,1]}\times\\
&(s_3,s_2)(s_2,s_1)T_{[s_3-1,s_1]}'\equiv(t_3,t_2)(t,s)V_{[t_2-1,1]}(t_2,t_1)W_{[t_3-1,t_1]}'R_{[s_2-1,1]}\times\\
&(s_2,s_1)(s_3,s_1)T_{[s_3-1,s_1]}'\equiv(t_3,t_2)(t,s)V_{[t_2-1,1]}(t_2,t_1)W_{[t_3-1,t_1]}'R_{[s_2-1,1]}\times\\
&(s_2,s_1)T_{[s_3-1,s_1]}(s_3,s_1).
\end{align*}

Now we consider compositions of intersection of (\ref{E6}) with
(\ref{E6}-\ref{E9}):

\begin{align*}
&(\ref{E6})\wedge(\ref{E6})(t_3,s)V_{[t_2-1,1]}(t_2,t_1)W_{[t_3-1,t_1]}(t_3,t_1)\times\\
&R_{[s_2-1,1]}(s_2,s_1)T_{[t_3-1,s_1]}(t_3,s_1), t_2>s, t_3>s_2,
s_2>t_1,\\
&(\ref{E7})\wedge(\ref{E6})(2,1)V_{2[2,1]}(3,1)\dots V_{n-1[n-1,1]}(n,1)\times\\
&R_{[t_2-1,1]}(t_2,t_1)T_{[n-1,t_1]}(n,t_1),\\
&(\ref{E6})\wedge(\ref{E8})(t_3,s)V_{[t_2-1,1]}(t_2,t_1)W_{[t_3-1,t_1]}(t_3,t_1)\delta^{\pm1},
t_2>s.
 \end{align*}

We check one:
\begin{align*}
&(\ref{E6})\wedge(\ref{E8})w=(n,s)V_{[t_2-1,1]}(t_2,t_1)W_{[n-1,t_1]}(n,t_1)\delta,
n=t_3, t_2>s,\\
&u-v=(t_2,s)(n,s)V_{[t_2-1,1]}(t_2,t_1)W_{[n-1,t_1]}'\delta -\\
&(n,s)V_{[t_2-1,1]}(t_2,t_1)W_{[n-1,t_1]}\delta(t_1+1,1),\\
&u\equiv\delta(t_2+1,s+1)(s+1,1)V_{[t_2,2]}(t_2+1,t_1+1)W_{[n,t_1+1]}^\star,\\
&\text { where } W_{[n,t_1+1]}(t_1+1,1)=(t_1+1,1)W_{[n,t_1+1]}^\star,\\
&v\equiv\delta(s+1,1)V_{[t_2,2]}(t_2+1,t_1+1)W_{[n,t_1+1]}(t_1+1,1)\equiv\\
&\delta(s+1,1)V_{[t_2,2]}(t_2+1,t_1+1)(t_1+1,1)W_{[n,t_1+1]}^\star\equiv\\
&\delta(s+1,1)V_{[t_2,2]}(t_2+1,1)(t_2+1,t_1+1)W_{[n,t_1+1]}^\star\equiv\\
&\delta(s+1,1)(t_2+1,1)V_{[t_2,2]}(t_2+1,t_1+1)W_{[n,t_1+1]}^\star\equiv\\
&\delta(t_2+1,s+1)(s+1,1)V_{[t_2,2]}(t_2+1,t_1+1)W_{[n,t_1+1]}^\star.
\end{align*}

To finish with the case of compositions of intersection, we only
need to consider
\begin{align*}
&(\ref{E7})\wedge(\ref{E8})w=(2,1)V_{2[2,1]}(3,1)\dots
V_{n-1[n-1,1]}(n,1)\delta^{\pm1}.
\end{align*}

We check one composition:
\begin{align*}
&(\ref{E7})\wedge(\ref{E8})w=(2,1)V_{2[2,1]}(3,1)\dots
V_{n-1[n-1,1]}(n,1)\delta,\\
&u-v=\delta V_{2[2,1]}'\dots V_{n-1[n-1,1]}'\delta-
(2,1)V_{2[2,1]}(3,1)\dots
V_{n-1[n-1,1]}\delta(2,1),\\
&u\equiv \delta^2V_{2[2,1]}'^+\dots V_{n-1[n-1,1]}'^+, \text { where
}R^+=R|_{(x,y)\mapsto(x+1,y+1)},\\
&v\equiv\delta(3,2)V_{2[2,1]}^+(4,2)\dots(n,2)
V_{n-1[n-1,1]}^+(2,1)\equiv\\
&\delta(2,1)(3,1)V_{2[2,1]}^{+\star}(4,2)\dots(n,1)
V_{n-1[n-1,1]}^{+\star}\equiv\\
&\delta^2(V_{2[2,1]}^{+\star })'\dots (V_{n-1[n-1,1]}^{+\star})'.
\end{align*}

It is easy to see that
$$
V_{[t-1,1]}'^+=(V_{[t-1,1]}^{+\star})', 2\leq t\leq n-1,
$$
identically, see also Lemma \ref{L00}. Indeed,
\begin{align*}
&V_{[t-1,1]}'^+=(V|_{(x,y)\mapsto(x,y),y\neq 1,
(x,1)\mapsto(t,x)})^+=V|_{(x,y)\mapsto(x+1,y+1),y\neq 1,
(x,1)\mapsto(t+1,x+1)},\\
&(V_{[t-1,1]}^{+\star})'=(((V|_{(x,y)\mapsto(x+1,y+1),y\neq 1,
(x,1)\mapsto
(x+1,2)})_{[t,2]})^\star)'=\\
&(V|_{(x,y)\mapsto(x+1,y+1),y\neq 1, (x,1)\mapsto
(x+1,1)})_{[t,2]}'=V|_{(x,y)\mapsto(x+1,y+1),y\neq 1, (x,1)\mapsto
(t+1,x+1)}.
\end{align*}

Now we proceed with compositions of including of relations
(\ref{E1})-(\ref{E9}). First of all, if a word
$(k_1,l_1)(i_1,j_1),k_1>l_1>i_1>j_1$ is a subword of words $V,
W,...$ of left parts of (\ref{E2})-(\ref{E7}), then  the triviality
of the corresponding compositions is clear. Let us list other
ambiguities,  using the same notations for (\ref{E1})-(\ref{E9}) as
before:
\begin{align*}
&(\ref{E2})\vee(\ref{E1})(k,l)(i_1,j_1)V_{[j-1]}(i,j), l>i_1>j_1,\\
&(\ref{E4})\vee(\ref{E1})(t_3,t_1)(i,j)V_{[t_2-1]}(t_3,t_2), t_1>i>j,\\
&(\ref{E5})\vee(\ref{E1})(t,s)(i,j)V_{[t_2-1,1]}(t_2,t_1)W_{[t_3-1,t_1]}(t_3,t_1),s>i>j,\\
&(\ref{E6})\vee(\ref{E1})(t_3,s)(i,j)V_{[t_2-1,1]}(t_2,t_1)W_{[t_3-1,t_1]}(t_3,t_1),s>i>j.\\
\end{align*}
All these cases are clear.

For (\ref{E2})-(\ref{E9}):
\begin{align*}
&(\ref{E2})\vee(\ref{E2})(k,l)V_{[j_1-1,1]}(i_1,j_1)W_{[j-1]}(i,j),j>i_1>j_1>l,\\
&(\ref{E4})\vee(\ref{E2})(t_3,t_1)V_{[j-1,1]}(i,j)W_{[t_2-1]}(t_3,t_2),t_2>i>j>t_1,\\
&(\ref{E5})\vee(\ref{E2})(t,s)V_{[j-1,1]}(i,j)W_{[t_2-1,1]}(t_2,t_1)R_{[t_3-1,t_1]}(t_3,t_1),t_2>i>j>s,\\
&(\ref{E5})\vee(\ref{E2})(t,s)V_{[t_2-1,1]}(t_2,t_1)W_{[j-1,t_1]}(i,j)R_{[t_3-1,t_1]}(t_3,t_1),
t_2>i>j>t_1,\\
&(\ref{E6})\vee(\ref{E2})(t_3,s)V_{[j-1,1]}(i,j)W_{[t_2-1,1]}(t_2,t_1)R_{[t_3-1,t_1]}(t_3,t_1),t_2>i>j>s,\\
&(\ref{E6})\vee(\ref{E2})(t_3,s)V_{[t_2-1,1]}(t_2,t_1)W_{[j-1,t_1]}(i,j)R_{[t_3-1,t_1]}(t_3,t_1),t_2>i>j>t_1 ,\\
&(\ref{E7})\vee(\ref{E2})(2,1)V_{2[2,1]}(3,1)\dots(k,1)W_{[j-1,1]}(i,j)R_{k[k,1]}\dots
V_{n-1[n-1,1]}(n,1),\\
&k>i>j>1.
\end{align*}

Now let us work out the second composition of (\ref{E6})-(\ref{E2}):
\begin{align*}
&(\ref{E6})\vee(\ref{E2})w=(t_3,s)V_{[t_2-1,1]}(t_2,t_1)W_{[j-1,t_1]}(i,j)R_{[t_3-1,t_1]}(t_3,t_1),t_2>i>j>t_1 ,\\
&u-v=(t_2,s)(t_3,s)V_{[t_2-1,1]}(t_2,t_1)W_{[j-1,t_1]}'(i,j)R_{[t_3-1,t_1]}'-\\
&(t_3,s)V_{[t_2-1,1]}(i,j)(t_2,t_1)W_{[j-1,t_1]}R_{[t_3-1,t_1]}(t_3,t_1),\\
&u\equiv(t_2,s)(t_3,s)V_{[t_2-1,1]}(t_2,t_1)(i,j)W_{[j-1,t_1]}'R_{[t_3-1,t_1]}'\\
(&\text{for } W_{[j-1,t_1]}'=(W_{[j-1,t_1]}')_{[t_3,t_1+1]} \text{
and we can assume that } j>t_1+1, \\
&\text {otherwise }W_{[j-1,t_1]}=1)\equiv(t_2,s)(t_3,s)V_{[t_2-1,1]}(i,j)(t_2,t_1)W_{[j-1,t_1]}'R_{[t_3-1,t_1]}',\\
&v\equiv(t_2,s)(t_3,s)V_{[t_2-1,1]}(i,j)(t_2,t_1)W_{[j-1,t_1]}'R_{[t_3-1,t_1]}'.\\
\end{align*}

For (\ref{E3})-(\ref{E9}):
\begin{align*}
&(\ref{E4})\vee(\ref{E3})(t_3,t_1)(t_1,t_1')V_{[t_2-1]}(t_3,t_2), t_1>t_1',\\
&(\ref{E5})\vee(\ref{E3})(t,s)(s,s')V_{[t_2-1,1]}(t_2,t_1)W_{[t_3-1,t_1]}(t_3,t_1),s>s',\\
&(\ref{E6})\vee(\ref{E3})(t_3,s)(s,s')V_{[t_2-1,1]}(t_2,t_1)W_{[t_3-1,t_1]}(t_3,t_1),s>s'.\\
\end{align*}

For (\ref{E4})-(\ref{E9}):
\begin{align*}
&(\ref{E5})\vee(\ref{E4})(t,s)V_{[t_2-1,1]}(t_2,t_1)W_{[t_2'-1,t_1]}(t_2,t_2')R_{[t_3-1,t_1]}(t_3,t_1),t_2>t_2'>t_1,\\
&(\ref{E6})\vee(\ref{E4})(t_3,s)V_{[t_2-1,1]}(t_2,t_1)W_{[t_2'-1,t_1]}(t_2,t_2')R_{[t_3-1,t_1]}(t_3,t_1),t_2>t_2'>t_1.
\end{align*}

Let us work out the first composition:

\begin{align*}
&(\ref{E5})\vee(\ref{E4})w=(t,s)V_{[t_2-1,1]}(t_2,t_1)W_{[t_2'-1,t_1]}(t_2,t_2')R_{[t_3-1,t_1]}(t_3,t_1),t_2>t_2'>t_1,\\
&u-v=(t_3,t_2)(t,s)V_{[t_2-1,1]}(t_2,t_1)W_{[t_2'-1,t_1]}'(t_2,t_2')R_{[t_3-1,t_1]}'-\\
&(t,s)V_{[t_2-1,1]}(t_2',t_1)(t_2,t_1)W_{[t_2'-1,t_1]}R_{[t_3-1,t_1]}(t_3,t_1),\\
&u\equiv(t_3,t_2)(t,s)V_{[t_2-1,1]}(t_2,t_1)(t_2,t_2')W_{[t_2'-1,t_1]}'R_{[t_3-1,t_1]}'\\
&(\text{for }W_{[t_2'-1,t_1]}'=(W_{[t_2'-1,t_1]}')_{[t_3,t_1+1]}
\text{ and } t_2'>t_1+1, \text{ otherwise }W_{[t_2'-1,t_1]}=1)\\
&\equiv(t_3,t_2)(t,s)V_{[t_2-1,1]}(t_2',t_1)(t_2,t_1)W_{[t_2'-1,t_1]}'R_{[t_3-1,t_1]}',\\
&v\equiv(t_3,t_2)(t,s)V_{[t_2-1,1]}(t_2',t_1)(t_2,t_1)W_{[t_2'-1,t_1]}'R_{[t_3-1,t_1]}'.\\
\end{align*}

Theorem 2.7 is proved.

\section{Corollaries}

Let $S\subset k\langle X \rangle$. A word $u$ is called
$S$-irreducible if $u\neq a\overline{s}b$, where $s\in S, a,b \in
X^*$. Let $Irr(S)$ be the set of all $S$-irreducible words.

Recall  Shirshov's Composition lemma (\cite{Sh62}, \cite{Bo72},
\cite{Bo76}) :

{\it Let $S$ be a \GS\ set in $k\langle X \rangle$. If $f\in
ideal(S)$, then $\overline{f}=a\overline{s}b, s\in S$. The converse
is also true.}

The Main corollary to this lemma is the following statement
(\cite{Sh62}, \cite{Bo72}, \cite{Bo76}):

{\it A subset $S \subset k\langle X \rangle $ is a \GS\ set iff the
set $Irr(S)$ is a linear basis for the algebra $k\langle X \rangle /
ideal(S)= k\langle X| S \rangle$ generated by $X$ with defining
relations $S$.}

Let $G=sgp\langle X|S \rangle$ be the semigroup generated by $X$
with defining relations $S$. Then $S$ is called a \GS\ basis of $G$
if $S$ is a \GS\ basis of the semigroup algebra $k(G)$, i.e., $S$ is
a \GS\ set in $ k\langle X \rangle$. It follows from the Main
corollary to the Composition lemma that in this case any word $u$ in
$X$ is equal in $G$ to a unique $S$-irreducible word $C(u)$, called
the normal (canonical) form of $u$.

Now let $S$ be the set of relations (\ref{E1})-(\ref{E9}), and let
$C(u)$ be a normal form of a word $u\in B_{n}$. Then $C(u)$ has a
form
$$
C(u)=\delta^k A,
$$
where $k\in \mathbb{Z}$, and  $A$ a positive $S$-irreducible word in
$a_(i,j)$'s, $C(A)=A$. Then   $A\neq \delta A_1$ for every positive
word $A_1$, otherwise $A=\delta C(A_1)$ identically, that is
impossible.

As a result, we have the following

\begin{cor} The $S$-irreducible normal form of a word of
$B_{n}$ in Birman-Ko-Lee-Garside generators coincides with the
Birman-Ko-Lee normal form of the word.
\end{cor}

\begin{proof}
Recall (\cite{BKL98}) that the Birman-Ko-Lee  normal form $G(u)$ of
$u\in B_{n}$ is
$$
G(u)=\delta^k A,
$$
where $u=G(u)$ in $B_{n}$, $k\in \mathbb{Z}$, and $A$ is a positive
word in $a_(i,j)$, $A\neq \Delta A_1$ for every positive word $A_1$,
and $A$ is the minimal word  with these properties. We have proved
that the $S$-irreducible normal form $C(u)$ has these properties.
\end{proof}

\begin{cor}[\cite{BKL98}] The semigroup of
positive braids $B_{n}^+$ in Birman-Ko-Lee generators can be
embedded into a group.
\end{cor}

\begin{proof}  It follows immediately from Theorem 2.7 that $B_{n}^+\subset B_{n}$.
\end{proof}

\end{document}